\documentclass[a4paper,12pt]{article}

\usepackage[intlimits]{amsmath}
\usepackage{amssymb}
\usepackage{amsthm}
\usepackage{mathrsfs}
\usepackage{cite}
\usepackage{color}
\usepackage{bm}

\newcommand{\la}{\langle}
\newcommand{\ra}{\rangle}
\newcommand{\pr}{\partial}
\newcommand{\dom}{\Omega}

\newcommand{\N}{{\mathcal{N}}}

\newcommand{\R}{\mathbb{R}}

\newcommand{\Ord}{\mathscr{O}}

\newcommand{\dd}{\,\text{d}}
\newcommand{\supp}{\mbox{supp\;}}

\newtheorem{thm}{Theorem}
\newtheorem{cor}{Corollary}
\newtheorem{dfn}{Definition}

\title{An inverse boundary value problem for the inhomogeneous porous medium equation}
\author{C\u{a}t\u{a}lin I. C\^{a}rstea\thanks{Department of Applied Mathematics, National Yang Ming Chiao Tung University, Hsinchu, Taiwan, R.O.C.; \mbox{email: catalin.carstea@gmail.com}} \and Tuhin Ghosh\thanks{Department of Mathematics, Universit\"at Bielefeld, 33615 Bielefeld, Germany; \mbox{email: tghosh@math.uni-bielefeld.de}}\and  Gen Nakamura\thanks{Department of Mathematics, Hokkaido University, Sapporo  060-0808, Japan; \mbox{email: gnaka@math.sci.hokudai.ac.jp}}}
\date{}

\begin{document}
\maketitle

\begin{abstract}
In this paper we establish uniqueness in the inverse boundary value problem for determining the two coefficients in the inhomogeneous porous medium equation $\epsilon\pr_tu-\nabla\cdot(\gamma\nabla u^m)=0$, with $m>1$, in dimension 3 or higher, which is a degenerate parabolic type quasilinear PDE. Our approach relies on using a Laplace transform to turn the original equation into a coupled family of nonlinear elliptic equations, indexed by the frequency parameter ($1/h$ in our definition) of the transform. A careful analysis of the asymptotic expansion in powers of $h$, as $h\to\infty$, of the solutions to the transformed equation, with special boundary data, allows us to obtain sufficient information to deduce the uniqueness result.
\end{abstract}

\section{Introduction}

We would like to begin by recalling the classical inverse boundary value problem of Calder\'on, introduced in \cite{Ca}. Suppose $\dom\subset\R^n$, with $n\geq 3$, is a a bounded, smooth domain and $\gamma\in C^\infty(\dom)$ is a bounded positive function with a positive lower bound. If we consider the (forward) boundary value problem
\begin{equation}\label{eq-C}
\left\{
\begin{array}{l} \nabla\cdot(\gamma(x)\nabla u(x))=0,\quad x\in\dom\\[5pt]u|_{\pr\dom}=g\in C^\infty(\pr\dom),\end{array}
\right.
\end{equation}
then we can introduce the Dirichlet-to-Neumann map (DN map)
\begin{equation}
\Lambda_\gamma(g)=\left.\gamma \pr_\nu u\right|_{\pr\dom}=\left.\gamma \nu\cdot\nabla u\right|_{\pr\dom},
\end{equation}
where $u$ is the solution of \eqref{eq-C}, and $\nu$ is the outer pointing unit normal on $\pr\dom$. Calder\'on's problem refers to the question of determining $\gamma$ from a knowledge of $\Lambda_\gamma$. The important subproblem of uniqueness is the question of showing that the correspondence $\gamma\to\Lambda_\gamma$ is, in fact, injective. This question has been answered affirmatively in \cite{SU}. Note that, by the linearity of the DN map, if $\Lambda_\gamma(g)$ is given for only non-negative Dirichlet data $g$, then the whole $\Lambda_\gamma$ operator, defined on complex valued functions, can be recovered. We will use this observation later in the paper.

The subject of this paper is the following boundary value problem for the (inhomogeneous) {\it porous medium equation} (PME)
\begin{equation}\label{eq}
\left\{\begin{array}{l}\epsilon(x)\pr_t u(t,x) -\nabla\cdot(\gamma(x)\nabla u^m(t,x))=0,\quad (t,x)\in (0,\infty)\times\dom,\\[5pt] u(0,x)=0,\quad u|_{[0,\infty)\times\pr\dom}=\phi,\quad u\geq0.\end{array}\right.
\end{equation}
Here $m>1$, and the coefficients $\epsilon$ and $\gamma$ are bounded positive functions, which also have strictly positive lower bounds. For convenience we will assume $\epsilon,\gamma\in C^\infty(\overline{\dom})$.
Note that consistency requires that $\phi(0,x)=0$, $\phi\geq0$. 

The Dirichlet-to-Neumann map associated to \eqref{eq} (informally) should be
\begin{equation}
\Lambda_{\epsilon,\gamma}^{PM} (\phi)=\gamma\pr_\nu u^m|_{[0,\infty)\times\pr\dom}.
\end{equation}
In order to rigorously define this map, we must first discuss the existence of solutions to \eqref{eq}. We will first introduce a notion of weak solutions on arbitrary, finite time intervals. Then we will see how these lead to a notion of solutions to \eqref{eq} for infinite time.

As its name suggests, the PME can be seen as a model for the flow of a gas through a porous medium. The function $u(t,x)$ being the {\it density} of the gas at time $t$ and position $x$. If $p$ is the {\it pressure}, then {\it Darcy's law} states that the velocity of the gas is $-\frac{k}{\mu}\nabla p$, where the parameter $k$ is the {\it permeability} of the medium and $\mu$ is the {\it dynamic viscosity} of the gas. The pressure is related to the density via a {\it state equation} of the form $p=p_0u^\sigma$, where $p_0$ is a constant {\it reference pressure} and $\sigma\geq1$ is the {\it polytropic exponent}. The value of this exponent depends on the gas but also on other factors, such as if the process being described is isothermal (then $\sigma=1$), adiabatic (then $\sigma>1$), etc. Conservation of mass requires that the following {\it continuity equation} holds
\begin{equation}
\epsilon\pr_t u+\nabla\cdot\left(\left(-\frac{k}{\mu}\right)u\nabla p\right)=0,
\end{equation}
where $0<\epsilon<1$ is the porosity of the medium (i.e. the fraction of the volume of the material that the gas can fill). Setting $\gamma=-\frac{\sigma kp_0}{(1+\sigma)\mu}$, $m=\sigma+1$, and  noting the state equation, we see that we obtain exactly equation \eqref{eq}. The parameters $\epsilon$ and $\gamma$ then depend on the particular gas considered, and also on the properties of the medium. 

From the purely mathematical point of view, the PME is interesting as a particularly simple nonlinear modification of the classical heat equation, which nevertheless can exhibit very different properties from its linear counterpart. It can be classified as a degenerate parabolic type quasilinear equation. Though it retains some of the features of the heat equation, such as having a maximum principle and a comparison principle, solutions to the PME are not typically real analytic in time, and the unique continuation property does not hold. To illustrate this, consider the following solutions to the PME with constant coefficients, known as Barenblatt solutions:
\begin{equation}
u(t,x)=t^{-\alpha}\left(C-k|x|^2t^{-2\beta}\right)_+^{\frac{1}{m-1}}
\end{equation}
($\alpha$, $\beta$, $k$ depend on $m$, $n$; $C$ is arbitrary.) These solutions vanish on a set with non-empty interior and the interface between the vanishing and non-vanishing regions can be said to move at a finite speed.

From the point of view of the study of inverse problems, the degenerate parabolic nature of the PME makes it a particularly interesting case to study, because known results of any kind for degenerate parabolic or elliptic PDE are still very few. We feel it is particularly interesting to see in what way the kinds of results that can be obtained may be constrained by finite speed of propagation behavior of solutions.

A good reference for the mathematical theory of the PME is the monograph \cite{V}. Though the particular case of the inhomogeneous equation is not treated there, most of what follows regarding the existence of solutions to the boundary value problem is closely based on the treatment of the homogeneous case in \cite{V}.

Before giving the specific definition of weak solutions which we will be working with, we must introduce a few notations. For $0<T\leq\infty$, let $Q_T:=(0,T)\times\Omega$ and $S_T:=(0,T)\times\partial\Omega$. Let 
\begin{equation}
C_{t}(Q_T)=\left\{\varphi\in C^\infty(\overline{Q_T}): \supp\varphi\cap \left[\{T\}\times\dom \right]=\emptyset\right\}
\end{equation} 
and let $H^1_{t}(Q_T)$ be the completion of this space in $H^1(Q_T)$. We will also denote by $H^{\frac{1}{2}}_{t}(S_T)$ the subspace of $H^{\frac{1}{2}}(S_T)$ that consists of traces of $H^1_{t}(Q_T)$ functions.
Similarly, we also introduce
\begin{equation}
C_{\diamond}(Q_T)=\left\{\varphi\in C^\infty(\overline{Q_T}): \supp\varphi\cap \left[S_T\cup\{T\}\times\dom \right]=\emptyset\right\},
\end{equation}
and its completion $H^1_{\diamond}(Q_T)$ in $H^1(Q_T)$.

Since we will be working with weak solutions, it is useful to recall a few things about boundary traces. We will denote by $\tau_\dom$ the boundary trace operator for the domain $\dom$. It is bounded between $H^1(\dom)$ and $H^{\frac{1}{2}}(\pr\dom)$. 

We will denote by $\tau_{Q_T}$ the boundary trace operator for $Q_T$ to the side boundary $S_T$. It is not hard to check that $\tau_{Q_T}$ is bounded from $L^2((0,T);H^1(\dom))$ to $L^2((0,T);H^{\frac{1}{2}}(\pr\dom))$. 

For our convenience later in the paper, we would like to consider the equation \eqref{eq} with an added source term
\begin{equation}\label{eqs}
\left\{\begin{array}{l}\epsilon(x)\pr_t u(t,x) -\nabla\cdot(\gamma(x)\nabla u^m(t,x))=f(t,x),\quad (t,x)\in Q_\infty,\\[5pt] u(0,x)=0,\quad u|_{[0,\infty)\times\pr\dom}=\phi(t,x),\quad u\geq0.\end{array}\right.
\end{equation}

\begin{dfn}\label{def-sol}
Suppose $0<T<\infty$. We say that $u\in L^\infty(Q_T)$, $u\geq0$, is a weak solution of \eqref{eqs} in $Q_T$ if it
satisfies the following conditions:
\begin{enumerate}
\item $\nabla (u^m)$ exists in the sense of distributions and $\nabla (u^m)\in L^2(Q_T)$;
\item for any test function $\varphi\in C_{\diamond}(Q_T)$ (or, equivalently, any $\varphi\in H^1_{\diamond}(Q_T)$)
\begin{equation} \int_{Q_T}\gamma\nabla\varphi\cdot\nabla (u^m)\dd t\dd x-\int_{Q_T}\epsilon\partial_t\varphi\, u\dd t\dd x =\int_{Q_T}\varphi f\dd t\dd x;
\end{equation}
\item $\tau_{Q_T}(u^m)=\phi^m$.
\end{enumerate}
\end{dfn}

Let the H\"older spaces $\mathcal{H}^{\frac{l}{2},l}(Q_T)$, $l\in\R_+\setminus\mathbb{N}$ be as defined in \cite[Chapter I, page 7]{La}. In section \ref{forward} we will prove the following result.
\begin{thm}\label{thm-forward}
For any $0<T<\infty$, if $\phi\geq0$, $\phi\in C(\overline{S_T})\cap \mathcal{H}^{1+\frac{\beta}{2},2+\beta}(S_T)$,  $\beta\in(0,1)$,  $\phi|_{t=0}=0$, and if $f\in L^\infty(Q_T)$, $f\geq0$, then there exists a unique weak solution $u$ of \eqref{eqs} in $Q_T$. This solution satisfies the energy estimate
\begin{equation}
||u^m||_{L^2((0,T);H^1(\dom))}\leq C(1+T)^{\frac{1}{2}}\left( ||\phi^m||_{W^{1,\infty}(S_T)}+||\phi||_{L^\infty(S_T)} +||f||_{L^\infty(Q_T)}\right),
\end{equation}
with a constant $C>0$ that depends on $\dom$, $m$, and the upper and lower bounds of $\epsilon$ and $\gamma$.
The solution $u$ satisfies the maximum principle
\begin{equation}
0\leq \sup_{Q_T}u \leq \sup_{S_T} \phi.
\end{equation}
If $\phi_1\leq\phi_2$, $f_1\leq f_2$ are as above and give rise to weak solutions $u_1$ and $u_2$, then
\begin{equation}
u_1\leq u_2.
\end{equation}
\end{thm}

Since the length of the time interval in the above theorem is arbitrary, it follows that we also have infinite time weak solutions, in the sense specified below.

\begin{cor}\label{cor-forward}
If  $\phi\geq0$, $\phi\in C(\overline{S_\infty})\cap \mathcal{H}^{1+\frac{\beta}{2},2+\beta}_{loc}(S_\infty)$,  $\beta\in(0,1)$,  $\phi|_{t=0}=0$, then there exists $u:Q_\infty\to[0,\infty)$ such that, for any $0<T<\infty$, $u$ is a weak solution of \eqref{eq} in $Q_T$ and
\begin{equation}
0\leq \sup_{Q_T}u \leq \sup_{S_T} \phi.
\end{equation}
\end{cor}

This allows us to rigorously define the Dirichlet-to-Neumann map for all $\phi$ satisfying the conditions in Corollary \ref{cor-forward}. If $u$ is the corresponding weak solution and $\psi\in H^1_{t}(Q_T)$, then we can define the Neumann data $\gamma\pr_\nu u^m|_{S_T}$ for each $0<T\leq\infty$ by 
\begin{equation}
\la \gamma\pr_\nu u^m|_{S_T}, \psi|_{S_T}\ra=\int_{Q_T}\left(\gamma \nabla\psi\cdot\nabla(u^m)-\epsilon\pr_t\psi u \right)\dd t\dd x.
\end{equation}
We see then that $\Lambda_{\epsilon,\gamma}^{PM} (\phi)\in \left(H^{\frac{1}{2}}_t(S_T)\right)'$.

Suppose now that we have two sets of coefficients, $\epsilon^{(i)}$, $\gamma^{(i)}$ and $\epsilon^{(ii)}$, $\gamma^{(ii)}$ that are as above. Let $M>0$. The main result of our paper is the following.
\begin{thm}\label{thm-inverse}
If $\Lambda_{\epsilon^{(i)},\gamma^{(i)}}^{PM} (\phi)=\Lambda_{\epsilon^{(ii)},\gamma^{(ii)}}^{PM} (\phi)$, for all  $0\leq \phi\leq M$, $\phi\in C(\overline{S_\infty})\cap \mathcal{H}^{1+\frac{\beta}{2},2+\beta}(S_\infty)$,  $\beta\in(0,1)$,  $\phi|_{t=0}=0$, then $\gamma^{(i)}=\gamma^{(ii)}$ and $\epsilon^{(i)}=\epsilon^{(ii)}$.
\end{thm}

Note that if $u(t,x)$ is a solution to the porous medium equation, and $\lambda>0$, then $u_\lambda(t,x)=\lambda^{\frac{1}{m-1}}u(\lambda t, x)$ is also a solution to the same equation. It follows that it is sufficient to prove Theorem \ref{thm-inverse} for one value of $M$, such as $M=2$.

We will prove Theorem \ref{thm-inverse} in section \ref{inverse}. Our method involves first going from an equation for the function $u$ to an equivalent problem for the function $v=u^m$. We then use a Laplace transform in time to introduce functions
\begin{equation}
V(h,x)=\int_0^\infty e^{-\frac{t}{h}} v(t,x)\dd t.
\end{equation}
These satisfy a family of equations, indexed by $h$, of the form
\begin{equation}
\nabla\cdot(\gamma(x)\nabla V(h,x))=\mathcal{N}(h,x),
\end{equation}
with the right hand side also being constructed from the solution $u$. In principle this is an infinite, coupled system of equations. One key technical part of the argument is to show that $V(h,x)$ satisfies decoupled elliptic inequalities of the form
\begin{equation}
0\leq \nabla\cdot(\gamma(x)\nabla V(h,x))\leq h^{-\frac{1}{m}}\epsilon(x)\left(V(h,x)\right)^{\frac{1}{m}}.
\end{equation}
If we choose boundary data of the form $\phi^m(t,x)=\chi(t)^m g(x)$, we are then able to derive the first two terms in the asymptotic expansion of $V(h,x)$ in powers of $h$, as $h\to\infty$. The corresponding expansion of the Dirichlet-to-Neumann map leads to separate integral identities that the differences $\gamma^{(i)}-\gamma^{(ii)}$ and $\epsilon^{(i)}-\epsilon^{(ii)}$ must satisfy. These quickly lead to the claimed uniqueness result.

The result presented in this paper fits into a series of other works on inverse boundary value problems for nonlinear equations. In fact, our method  resembles the ``second linearization'' approach that has been used in the past for (non-degenerate) elliptic and parabolic semilinear and quasilinear problems. We list here some of the important papers in this field. For semilinear equations, examples include \cite{CaKi}, \cite{ChKi}, \cite{FO}, \cite{I1}, \cite{IN}, \cite{IS}, \cite{KiUh}, \cite{KU}, \cite{KU2}, \cite{LLLS1}, \cite{LLLS2}, \cite{S2}. For quasilinear equations, not in divergence form, see \cite{I2}. For quasilinear equations in divergence form see \cite{CF1}, \cite{CF2}, \cite{CFKKU}, \cite{CK}, \cite{CNV}, \cite{EPS}, \cite{HS}, \cite{KN}, \cite{MU}, \cite{Sh}, \cite{S1}, \cite{S3}, \cite{SuU}, (also \cite{C} for quasilinear time-harmonic Maxwell systems). 

Results for degenerate equations, such as the PME, are quite few. Examples would include the following works for the weighted p-Laplace equation, \mbox{$\nabla\cdot(a(x)|\nabla u|^{p-2}\nabla u)=0$}, which is a degenerate elliptic quasilinear PDE: \cite{BKS}, \cite{B}, \cite{BHKS}, \cite{GKS}, \cite{BIK}, \cite{KW}. Note that uniqueness in the inverse boundary value problem for the coefficient $a$ in this equation, without additional assumptions such as monotonicity (as in \cite{GKS}), has not yet been shown. 

To the best of our knowledge, this paper has provided the first unconstrained uniqueness result for a degenerate parabolic or elliptic quasilinear PDE. After the first version of this work has been made available, a second paper \cite{CGU} on inverse boundary value problems for the PME, by some of the same authors, has been published. In that work uniqueness is obtained with an arbitrarily small time of observation $T$. There is however a crucial difference: while in Theorem \ref{thm-inverse} we only need to work with bounded boundary data, for the analogous result in \cite{CGU} arbitrarily large boundary data is required. We conjecture that finite speed of propagation is the reason for this difference: one may ensure that the solutions generated by boundary data probe the whole domain either by waiting an infinite amount of time, or (seeing that larger Barenblatt solutions propagate faster) by increasing the size of the data/solutions. 

Finally, we would like to comment briefly on the related inverse problem for the {\it fast diffusion equation}, i.e. the equation \eqref{eq} with $0<m<1$. One expects, since finite speed of propagation phenomena do not occur, that the inverse boundary value problem should be easier to solve in this case. Indeed, we can give here  a sketch of how the inverse boundary value problem for the fast diffusion equation may be reduced to a known result. As for the PME, one starts with the change of function $v(t,x)=u(t,x)^m$,  which puts the equation in the form 
\begin{equation}
\epsilon\pr_t v^{\frac{1}{m}}-\nabla\cdot(\gamma\nabla v)=0.
\end{equation}
If we look for solutions of the form $v(t,x)=t^\alpha w(x)$, a simple computation gives that we must have  $\alpha=\frac{m}{1-m}>0$ and
\begin{equation}
-\nabla\cdot(\gamma\nabla w)+\frac{\epsilon}{1-m} w^{\frac{1}{m}}=0,
\end{equation}
so the question is reduced to the invese boundary value problem for this last equation.  This uniquness for the coefficients in this problem has already (almost) been established  in \cite{LLST}.

\section{The forward problem}\label{forward}

In this section we give a proof to Theorem \ref{thm-forward}, which gives the existence of finite time weak solutions. The method used is very similar to what can be found in \cite{V}, the difference being that in that reference the coefficients are taken to be constant.

We begin by introducing a sequence $u_k$, $k=1,2,\ldots$, consisting of solutions to the boundary value problem
\begin{equation}\label{eq-k}
\left\{\begin{array}{l}\epsilon\pr_t u_k -\nabla\cdot(\gamma\nabla u^m_k)=f_k,\\[5pt] u_k(0,x)=\frac{1}{k},\quad u_k|_{[0,T)\times\pr\dom}=\phi+\frac{1}{k},\end{array}\right.
\end{equation}
where we choose functions $f_k\in C^\infty(Q_T)$ such that
\begin{equation}
f(t,x)\leq f_{k+1}(t,x)\leq f_k(t,x)\leq f(t,x)+\frac{1}{k}.
\end{equation}
Intuitively, we expect the solutions $u_k$ (if they exist) to satisfy the maximum principle
\begin{equation}\label{mp-k}
\frac{1}{k}\leq u_k(t,x)\leq\frac{1}{k}+\sup_{S_T}\phi+\mu t\sup_{Q_T}f_k,
\end{equation}
where $\mu=(\inf_{\dom}\epsilon)^{-1}$.
With this in mind, we redefine the $u_k$ to be the solutions of 
\begin{equation}\label{eq-k-2}
\left\{\begin{array}{l}\epsilon\pr_t u_k -\nabla\cdot\left(a_k(x,u_k)\nabla u_k)\right)=f_k,\\[5pt] u_k(0,x)=\frac{1}{k},\quad u_k|_{[0,T)\times\pr\dom}=\phi+\frac{1}{k},\end{array}\right.
\end{equation}
where the function $a_k(x,\lambda)$ is chosen so as to make \eqref{eq-k-2} a uniform quasilinear parabolic problem (i.e. it avoids the degeneracy at $\lambda=0,\infty$), but such that if 
\begin{equation}\frac{1}{k}\leq\lambda\leq \sup_{S_T}\phi+\mu T\sup_{Q_T} f+\frac{1+\mu T}{k},
\end{equation}
we have
\begin{equation}
a_k(x,\lambda)=m\gamma(x)\lambda^{m-1}.
\end{equation}
Outside this domain $a_k$ is extended smoothly in such a way that \eqref{eq-k-2} is non-degenerate.
In this case (see \cite[Chapter V, page 452, Theorem 6.1]{La}), the problem \eqref{eq-k-2} has a solution $u_k \in C^{1,2}(Q_T)\cap C(\overline{Q_T})$ which furthermore 
satisfies (see \cite[Theorem 9.7]{Li}) the maximum principle \eqref{mp-k}. Since \eqref{mp-k} is satisfied, it follows that $u_k$ is in fact also a solution to \eqref{eq-k}.

Note that $u_{k}$ is also a solution to $\epsilon\pr_t u_{k} -\nabla\cdot\left(a_{k+1}(x,u_k,\nabla u_k)\right)=f_k$. Therefore, by the comparison principle (see \cite[Theorem 9.7]{Li}), we have that
\begin{equation}
0\leq u_{k+1}\leq u_{k},\quad\forall k=1,2\ldots.
\end{equation}
We can then define the pointwise limit
\begin{equation}
u(t,x)=\lim_{k\to\infty} u_k(t,x).
\end{equation}
This function will turn out to be the weak solution to \eqref{eq} that we are searching for.

Since clearly $u\leq u_k$ for all $k$, by \eqref{mp-k} we have that 
\begin{equation}
0\leq u(t,x)\leq \sup_{S_T}\phi+\mu t\sup_{Q_T} f.
\end{equation}
By the monotone convergence theorem we have that
\begin{equation}
\lim_{k\to\infty}\int_{Q_T}(u_k-u)^2\dd t\dd x=0,
\end{equation}
i.e. that $u_k\to u $ in $L^2(Q_T)$. Similarly, it also holds that $u_k^m\to u^m $ in $L^2(Q_T)$.

We will denote by $\tilde\phi:Q_T\to [0,\infty)$ a smooth extension of the Dirichlet data, such that we still have $\tilde\phi(0,x)=0$ for all $x\in\dom$. We can choose this extension in such a way that
\begin{equation}
||\tilde\phi||_{W^{1,\infty}(Q_T)}\leq C||\phi||_{W^{1,\infty}(S_T)},
\end{equation}
with a constant $C>0$.

Let 
\begin{equation}
\eta_k=u_k^m-\left(\tilde\phi+\frac{1}{k}\right)^m,
\end{equation}
which is zero on $S_T$. We can multiply \eqref{eq-k} by $\eta_k$ and integrate over $Q_T$. Note that
\begin{multline}
-\int_{Q_T}\eta_k\nabla(\gamma\nabla u_k^m)\dd t\dd x\\[5pt]
=\int_{Q_T}\gamma|\nabla u_k^m|^2\dd t\dd x
-\int_{Q_T}\gamma\nabla \left(\tilde\phi+\frac{1}{k}\right)^m\cdot\nabla u_k^m\dd t\dd x,
\end{multline}
and that
\begin{multline}
\int_{Q_T}\epsilon\pr_t u_k\eta_k\dd t\dd x
=\int_{Q_T}\epsilon\pr_t u_k u_k^m\dd t\dd x
-\int_{Q_T}\epsilon\pr_t u_k\left(\tilde\phi+\frac{1}{k}\right)^m\dd t\dd x\\[5pt]
=\int_{\dom}\frac{\epsilon}{m+1}u_k^{m+1}(T)\dd x
-\int_{\dom}\epsilon u_k(T)\left(\tilde\phi(T)+\frac{1}{k}\right)^m\dd x\\[5pt]
+\int_{Q_T}\epsilon u_k\pr_t\left(\tilde\phi+\frac{1}{k}\right)^m\dd t\dd x.
\end{multline}
Putting these two together we have 
\begin{multline}
\int_{Q_T}\gamma|\nabla u_k^m|^2\dd t\dd x
+\int_{\dom}\frac{\epsilon}{m+1}u_k^{m+1}(T)\dd x\\[5pt]
=\int_{Q_T}\gamma\nabla \left(\tilde\phi+\frac{1}{k}\right)^m\cdot\nabla u_k^m\dd t\dd x
+\int_{\dom}\epsilon u_k(T)\left(\tilde\phi(T)+\frac{1}{k}\right)^m\dd x\\[5pt]
-\int_{Q_T}\epsilon u_k\pr_t\left(\tilde\phi+\frac{1}{k}\right)^m\dd t\dd x
+\int_{Q_T}f_k\left[u_k^m-\left(\tilde\phi+\frac{1}{k}\right)^m\right].
\end{multline}
By straightforward estimates we can conclude that there exists a constant $C>0$, depending on $m$ and the upper and lower bounds of $\epsilon$ and $\gamma$, such that
\begin{multline}
\int_{Q_T}|\nabla u_k^m|^2\dd t\dd x+\int_\dom u_k^{m+1}(T)\dd x
\leq C\Bigg[||u_k||_{L^2(Q_T)}^2 +||u_k^m||_{L^2(Q_T)}^2+||f_k||_{L^2(Q_T)}^2\\[5pt]
+\left\Vert\left(\tilde\phi+\frac{1}{k}\right)^m\right\Vert_{H^1(Q_T)}^2+\left\Vert\left(\tilde\phi+\frac{1}{k}\right)^m\right\Vert_{L^\frac{m+1}{m}(\dom)}
\Bigg].
\end{multline}
There is then a constant $C'>0$, independent of $k$ and $T$, such that
\begin{multline}
\int_{Q_T}|\nabla u_k^m|^2\dd t\dd x<C'(1+T)\Bigg(\left\Vert\left(\phi+\frac{1}{k}\right)^m\right\Vert_{C^{0,1}(S_T)}^{2}
+\left\Vert\left(\phi+\frac{1}{k}\right)^m\right\Vert_{C(S_T)}\\[5pt] +\left\Vert\phi+\frac{1}{k}\right\Vert_{C(S_T)}^2 +\left\Vert f+\frac{1}{k}\right\Vert_{L^\infty(Q_T)}^2 \Bigg).
\end{multline}
It follows that (a subsequence of) $\nabla u_k^m$ converges weakly in $L^2(Q_T)$ to a limit $U$. It is easy to see that $U=\nabla u^m$ in the sense of distributions.

Since we have that $\tau_{Q_T}(u_k^m)(t)\to \phi^m$ in $L^2((0,T); H^{\frac{1}{2}}(\pr\dom))$ by construction, and also that $\tau_{Q_T}(u_k^m)\rightharpoonup \tau_{Q_T}(u^m)$, it follows that
\begin{equation}
\tau_{Q_T}(u^m)=\phi^m.
\end{equation}
This shows that $u$ is a weak solution of \eqref{eq} in $Q_T$.

In order to prove uniqueness of weak solutions, we start by assuming $u_1$ and $u_2$ are both weak solutions of \eqref{eq} in $Q_T$. We have that for all test functions $\varphi$
\begin{equation}
\int_{Q_T}\gamma (\nabla u_1^m-\nabla u_2^m)\nabla\varphi\dd t\dd x-
\int_{Q_T}\epsilon (u_1-u_2)\pr_t\varphi\dd t\dd x=0.
\end{equation}
Let
\begin{equation}
\varphi(t,x)=\int_t^T(u_1^m(s,x)-u_2^m(s,x))\dd s.
\end{equation}
It is not hard to check that $\varphi\in H^1_\diamond(Q_T)$, so we can use it in the above identity. We obtain
\begin{equation}
\frac{1}{2}\int_\dom\gamma\left(\int_0^T(\nabla u_1^m-\nabla u_2^m)\dd t \right)^2\dd x
+\int_{Q_T}\epsilon(u_1-u_2)(u_1^m-u_2^m)\dd t\dd x=0.
\end{equation}
The first term is clearly positive and, by the monotonicity of the $z\to z^m$ function in $\R_+$, so is the second. It follows that $u_1=u_2$.

Regarding the energy inequality claimed in the statement of the Theorem, it is enough to note that 
\begin{multline}
||\nabla u^m||_{L^2(Q_T)}\leq \liminf_{k\to\infty} ||\nabla u_k^m||_{L^2(Q_T)}\\[5pt]
\leq C(1+T)^{\frac{1}{2}}\left(||\phi^m||_{W^{1,\infty}(S_T)}+||\phi||_{L^\infty(S_T)} + ||f||_{L^\infty(Q_T)}\right).
\end{multline}

Finally, suppose we have boundary data and sources $\phi_1\leq\phi_2$, $f_1\leq f_2$, with corresponding solutions $u_1$ and $u_2$. It is enough to observe that the approximating sequences defined above must satisfy
\begin{equation}
u_{1,k}\leq u_{2,k},
\end{equation}
and this property is preserved in the limit.
This concludes the proof.

\section{The inverse problem}\label{inverse}

In this section we give a proof of Theorem \ref{thm-inverse}. We do this by transforming the equation so that the principal term becomes linear. Then we use a Laplace transform in time to transform the equation into a (coupled) family of nonlinear elliptic equations. Below, when we write $\epsilon$ and $\gamma$, we mean any one of the pairs $\epsilon^{(i)}$, $\gamma^{(i)}$ or $\epsilon^{(ii)}$, $\gamma^{(ii)}$.

\subsection{The first reformulation of the problem}

We first perform the change of function $v=u^m$. The new function satisfies the equation
\begin{equation}\label{eq-v}
\left\{\begin{array}{l}\epsilon(x)\pr_t v(t,x)^{\frac{1}{m}} -\nabla\cdot(\gamma(x)\nabla v(t,x))=0,\\[5pt] v(0,x)=0,\quad v|_{[0,\infty)\times\pr\dom}=f(t,x),\quad v\geq0,\end{array}\right.
\end{equation}
where $f=\phi^m$. It is clear that the notions of weak solutions we have introduced for equation \eqref{eq} give corresponding notions of weak solutions to \eqref{eq-v}, whose existence is then guaranteed by Theorem \ref{thm-forward}. The associated Dirichlet-to-Neumann map is 
\begin{equation}
\Lambda^v_{\epsilon,\gamma} (f)=\gamma\pr_\nu v|_{[0,\infty)\times\pr\dom},
\end{equation}
and it contains the same information as $\Lambda_{\epsilon,\gamma}^{PM}$, when considering $f\in C(\overline{S_\infty})\cap C^\infty(S_\infty)$, $f\geq 0$, and such that $f(0,x)=0$ for all $x\in\dom$.

\subsection{Reformulation of the problem using the Laplace transform}

The next reformulation of the problem uses the Laplace transform to turn our equation into an elliptic one. From now on we will only consider boundary data of a particular form, namely $f(t,x)=\chi(t)^mg(x)$, with $g\in C(\pr\dom)$, and $\chi\in C^\infty([0,\infty))$ an increasing function such that $\chi(t)=0$ if $t\leq1/2$ and $\chi(t)=1$ if $t\geq1$. For $h>0$ define
\begin{equation}
V(h,x)=\int_0^\infty e^{-\frac{t}{h}} v(t,x)\dd t.
\end{equation}
Note that since (in the sense of distributions) we have
\begin{equation}
\nabla V(h,x)=\int_0^\infty e^{-\frac{t}{h}} \nabla v(t,x)\dd t,
\end{equation}
it then follows that 
\begin{multline}
||\nabla V(h,\cdot)||_{L^2(\dom)}\leq \int_0^\infty e^{-\frac{t}{h}}||\nabla v(t,\cdot)||_{L^2(\dom)}\dd t\\[5pt]
\leq \left( \int_0^\infty e^{-\frac{t}{h}}\dd t \right)^\frac{1}{2}\left(\;\int_{Q_\infty} e^{-\frac{t}{h}}|\nabla v(t,x)|^2\dd t\dd x\right)^\frac{1}{2}<\infty,
\end{multline}
where the last inequality holds since, by Theorem \ref{thm-forward},
\begin{equation}
\int_{Q_T}|\nabla v(t,x)|^2\dd t\dd x
\end{equation}
has an at most polynomial growth as $T\to\infty$, so the integral over $Q_\infty$ converges.
Therefore we have that $V(h,\cdot)\in H^1(\dom)$. 

If $\varphi\in C_0^\infty(\dom)$, and $\zeta\in C_0^\infty(\R)$ is such that $\zeta(t)=1$ when $|t|\leq 1$ and $\zeta(t)=0$ when $|t|\geq 2$, then
\begin{multline}
\left\la \nabla(\gamma\nabla V(h,\cdot)),\varphi\right\ra=-\int_\dom \gamma(x)\nabla V(h,x)\cdot\nabla\varphi(x)\dd x\\[5pt]
=-\lim_{T\to\infty}\int_{Q_\infty}\zeta(\frac{t}{T})e^{-\frac{t}{h}}\gamma(x)\nabla v(t,x)\cdot\nabla\varphi(x)\dd t\dd x\\[5pt]
=-\lim_{T\to\infty}\int_{Q_\infty}\epsilon(x)\left[\frac{1}{T}\zeta'(\frac{t}{T})e^{-\frac{t}{h}}-\frac{1}{h}\zeta(\frac{t}{T})e^{-\frac{t}{h}} \right]v^\frac{1}{m}(t,x) \varphi(x)\dd t\dd x\\[5pt]
= \int_\dom\varphi(x)\left[\frac{1}{h}\epsilon(x)\int_0^\infty e^{-\frac{t}{h}} v^\frac{1}{m}(t,x)\dd t \right]\dd x.
\end{multline}
It follows that 
\begin{equation}
\nabla\cdot(\gamma(x)\nabla V(h,x))=\mathcal{N}(h,x),
\end{equation}
where
\begin{equation}
\mathcal{N}(h,x)
=\frac{1}{h}\epsilon(x)\int_0^\infty e^{-\frac{t}{h}} v^{\frac{1}{m}}(t,x)\dd t.
\end{equation}
The value $V(h,x)$ takes on the boundary (in the sense of traces) is
\begin{equation}
X(h)g(x)=\left(\int_0^\infty e^{-\frac{t}{h}} \chi(t)^m \dd t\right)g(x).
\end{equation}
Here we note that
\begin{equation}
he^{-\frac{1}{h}}=\int_1^\infty e^{-\frac{t}{h}}\dd t\leq X(h)\leq \int_0^\infty e^{-\frac{t}{h}}\dd t=h,
\end{equation}
so $X(h)=\Ord(h)$, $X(h)^{-1}=\Ord(h^{-1})$ as $h\to\infty$.

We then conclude that $V$ satisfies the family of boundary value problems
\begin{equation}\label{Laplace-elliptic}
\left\{\begin{array}{l} \nabla\cdot(\gamma(x)\nabla V(h,x))=\mathcal{N}(h,x),\\[5pt] V|_{\pr\dom}=X(h)g,\quad V\geq0.\end{array}\right.
\end{equation}
We define the associated DN map as essentially (up to a multiplicative correction) the Laplace transform of $\Lambda^v_{\epsilon,\gamma}$, namely
\begin{equation}
\Lambda^h_{\epsilon,\gamma}(g)=X(h)^{-1}\gamma\pr_\nu V(h,\cdot)|_{\pr\dom}.
\end{equation}
We have added the $X(h)^{-1}$ factor for the sake of convenience in the computations that follow below.

Let $m'$ be such that $m'^{-1}+m^{-1}=1$. Then by H\"older's inequality
\begin{multline}\label{H-trick}
0\leq \mathcal{N}(h,x)=\frac{1}{h}\epsilon(x)\int_0^\infty e^{-\frac{t}{m'h}}\left(e^{-\frac{t}{h}} v(t,x)\right)^{\frac{1}{m}}\dd t\\[5pt]
\leq h^{-\frac{1}{m}}\epsilon(x)\left(V(h,x)\right)^{\frac{1}{m}},
\end{multline}
so
\begin{multline}
||\mathcal{N}(h,\cdot)||_{L^\infty(\dom)}\leq h^{-\frac{1}{m}} ||\epsilon||_{L^\infty(\dom)}||V(h,\cdot)||_{L^\infty(\dom)}^{\frac{1}{m}}\\[5pt]\leq C h^{-\frac{1}{m}}||V(h,\cdot)||_{W^{1,p}(\dom)}^{\frac{1}{m}},
\end{multline}
for a $p>n$.

By eliptic regularity estimates for \eqref{Laplace-elliptic} we have
\begin{multline}
||V(h,\cdot)||_{W^{2,p}(\dom)}\leq C\left(h||g||_{W^{2-\frac{1}{p},p}(\pr\dom)}+||\mathcal{N}(h,\cdot)||_{L^p(\dom)} \right)\\[5pt]
\leq C\left(h||g||_{W^{2-\frac{1}{p},p}(\pr\dom)}+ h^{-\frac{1}{m}}||V(h,\cdot)||_{W^{1,p}(\dom)}^{\frac{1}{m}} \right).
\end{multline}
Then
\begin{multline}
\max\left(||V(h,\cdot)||_{W^{2,p}(\dom)},1\right)\\[5pt]\leq 
C\left(h||g||_{W^{2-\frac{1}{p},p}(\pr\dom)}+1+ h^{-\frac{1}{m}}\left[\max\left(||V(h,\cdot)||_{W^{1,p}(\dom)},1\right)\right]^{\frac{1}{m}} \right)\\[5pt]
\leq C\left(h||g||_{W^{2-\frac{1}{p},p}(\pr\dom)}+1+ h^{-\frac{1}{m}}\max\left(||V(h,\cdot)||_{W^{1,p}(\dom)},1\right) \right).
\end{multline}
There clearly exists an $h_0>0$ such that if $h>h_0$ we have
\begin{equation}
||V(h,\cdot)||_{W^{2,p}(\dom)}\leq h C \left(||g||_{W^{2-\frac{1}{p},p}(\pr\dom)}+h^{-1}\right).
\end{equation}

\subsection{First order asymptotic expansions and the determination of $\gamma$}

We make the following Ansatz
\begin{equation}\label{ansatz1}
V(h,x)=X(h)V_0(x)+
R_1(h,x),
\end{equation}
where $V_0\in W^{2,p}(\dom)$ solves the equation
\begin{equation}
\left\{\begin{array}{l} \nabla\cdot(\gamma(x)\nabla V_0(x))=0,\\[5pt] V_0|_{\pr\dom}=g,\quad g\geq0.\end{array}\right.
\end{equation}
It follows that $R_1$ must satisfy
\begin{equation}
\left\{\begin{array}{l} \nabla\cdot(\gamma(x)\nabla R_1(h,x))=
\mathcal{N}(h,x),\\[5pt] R_1|_{\pr\dom}=0.\end{array}\right.
\end{equation}
We observe that $R_1$ is a subsolution for the elliptic operator $\nabla\cdot\gamma\nabla$, so we must have that $R_1\leq 0$.

By elliptic estimates we have that $R_1$ is bounded uniformly in $h>h_0$
\begin{multline}\label{R1-estimate}
||R_1(h,\cdot)||_{W^{2,p}(\dom)}\leq C 
||\mathcal{N}(h,\cdot)||_{L^p(\dom)}\\[5pt]
\leq C h^{-\frac{1}{m}}||V(h,\cdot)||_{W^{1,p}(\dom)}^{\frac{1}{m}}
\leq C \left(||g||_{W^{2-\frac{1}{p},p}(\pr\dom)}^{\frac{1}{m}}+h^{-\frac{1}{m}}\right).
\end{multline}
The asymptotic behavior of $\Lambda^h_{\epsilon,\gamma}(g)$ as $h\to\infty$ is then
\begin{equation}\label{DN-expansion}
\Lambda^h_{\epsilon,\gamma}(g)=\gamma\pr_\nu V_0(\cdot)|_{\pr\dom}+\Ord(h^{-1}).
\end{equation}

Now, since we must have $\Lambda^h_{\epsilon^{(i)},\gamma^{(ii)}}(g)=\Lambda^h_{\epsilon^{(i)},\gamma^{(ii)}}(g)$ for all $g\in C(\overline{\dom})\cap C^\infty(\dom)$, with $g\geq0$, it follows that corresponding orders in $h$ of their asymptotic expansions must also be equal. From the $h^0$ order term we then get that the DN maps for the original Calder\'on problem, with conductivities $\gamma^{(i)}$ and $\gamma^{(ii)}$, are equal. That is, we have $\Lambda_{\gamma^{(i)}}=\Lambda_{\gamma^{(ii)}}$. It follows then form \cite{SU} that $\gamma^{(i)}=\gamma^{(ii)}$.

\subsection{Asymptotic expansions up to the second order and the determination of $\epsilon$}

In order to also prove the equality of the $\epsilon^{(i)}$, $\epsilon^{(ii)}$ coefficients, we need to look deeper into the asymptotic expansion of the DN map. Noting that
\begin{equation}
X(h)V_0(x)=\int_0^\infty e^{-\frac{t}{h}}\chi(t)^mV_0(x)\dd t,
\end{equation}
and also that $V_0\geq0$, let
\begin{equation}
\mathcal{N}_0(h,x)=\frac{1}{h}\epsilon(x)\int_0^\infty e^{-\frac{t}{h}}\chi(t)V_0^{\frac{1}{m}}(x)\dd t=\epsilon(x) \hat\chi(h)V_0^{\frac{1}{m}}(x),
\end{equation}
where
\begin{equation}
\hat\chi(h)=\frac{1}{h}\int_0^\infty e^{-\frac{t}{h}}\chi(t)\dd t=\Ord(1).
\end{equation}
We refine the Ansatz \eqref{ansatz1} to 
\begin{equation}\label{ansatz2}
V(h,x)=X(h)V_0(x)+V_1^h(x)+R_2(h,x),
\end{equation}
where $V_1^h\in W^{2,p}(\dom)$ is the solution to
\begin{equation}
\left\{\begin{array}{l} \nabla\cdot(\gamma(x)\nabla V_1^h(x))=\mathcal{N}_0(h,x),\\[5pt] V_1^h|_{\pr\dom}=0.\end{array}\right.
\end{equation}
As $V_1^h$ is a subsolution for the elliptic operator $\nabla\cdot\gamma\nabla$, it follows that $V_1\leq 0$.


Returning to \eqref{ansatz2}, we see that $R_2$ must satisfy the equation
\begin{equation}
\left\{\begin{array}{l} \nabla\cdot(\gamma(x)\nabla R_2(h,x))=
\left(\mathcal{N}(h,x)-\mathcal{N}_0(h,x)\right),\\[5pt] R_2|_{\pr\dom}=0.\end{array}\right.
\end{equation}
Let
\begin{equation}
w_0(t,x)=\chi(t)^mV_0(x)
\end{equation}
and note that
\begin{equation}
\left\{\begin{array}{l}
\epsilon\pr_t w_0^{\frac{1}{m}}-\nabla\cdot(\gamma\nabla w_0)=\epsilon \chi'(t)V_0^{\frac{1}{m}}\geq0,\\[5pt]
w_0|_{t=0}=0, w_0|_{(0,\infty)\times\pr\dom}=v|_{(0,\infty)\times\pr\dom}.\end{array}\right.
\end{equation}
 By the last part of Theorem \ref{thm-forward}, it follows that $w_0\geq v$. This implies that
\begin{equation}
\mathcal{N}(h,x)\leq\mathcal{N}_0(h,x).
\end{equation}
 It follows that $R_2$ is a supersolution for the elliptic operator $\nabla\cdot\gamma\nabla$, so we have that $R_2\geq0$.

We would like to also obtain an estimate for the difference $\N_0-\N$. To this end, using the same H\"older inequality trick employed above in equation \eqref{H-trick}, we estimate
\begin{multline}
0\leq \N_0(h,x)-\N(h,x)=\frac{\epsilon(x)}{h}\int_0^\infty e^{-\frac{t}{h}}(w_0^{\frac{1}{m}}(t,x)-v^{\frac{1}{m}}(t,x))\dd t\\[5pt]
\leq \frac{\epsilon(x)}{h}\int_0^\infty e^{-\frac{t}{h}} (w_0(t,x)-v(t,x))^{\frac{1}{m}}\dd t\\[5pt]
\leq \epsilon(x) h^{-\frac{1}{m}}\left(X(h)V_0(x)-V(h,x)\right)^{\frac{1}{m}}\\[5pt]
=\epsilon(x) h^{-\frac{1}{m}}\left(-R_1(h,x)\right)^{\frac{1}{m}}.
\end{multline}
It follows from \eqref{R1-estimate} that 
\begin{equation}
||\N-\N_0||_{L^\infty(\dom)}=\Ord(h^{-\frac{1}{m}}),
\end{equation}
as $h\to\infty$.
By elliptic estimates
\begin{equation}
||R_2(h,\cdot)||_{W^{2,p}(\dom)}=\Ord(h^{-\frac{1}{m}}).
\end{equation}
This gives the following asymptotic expansion for $\Lambda^h_{\epsilon,\gamma}$ as $h\to\infty$
\begin{equation}\label{DN-expansion2}
\Lambda^h_{\epsilon,\gamma}(g)=\gamma\pr_\nu V_0(\cdot)|_{\pr\dom}+h^{-1}\gamma\pr_\nu V_1(\cdot)|_{\pr\dom}+
\Ord(h^{-1-\frac{1}{m}}).
\end{equation}
As before, we will use  the $h^{-1}$ order in the expansion in order to obtain the uniqueness result for $\epsilon$.

Let $W$ be any smooth solution to $\nabla\cdot(\gamma\nabla W)=0$. We have that
\begin{multline}
\la \gamma\pr_\nu V_1|_{\pr\dom},W|_{\pr\dom}\ra
=\int_\dom\gamma\nabla V_1\cdot\nabla W\dd x+
\int_\dom\epsilon V_0^{\frac{1}{m}}W\dd x\\[5pt]
=\int_\dom\epsilon V_0^{\frac{1}{m}}W\dd x.
\end{multline}
If $\Lambda^h_{\epsilon,\gamma}$ is known, then so are all the integral quantities of the form that appears on the right hand side of the above equation.

We can take $V_0$ to be of the form
\begin{equation}
V_0(x)=1+sH(x),\quad \nabla\cdot(\gamma\nabla H)=0,\;H>0.
\end{equation}
Then
\begin{equation}
m\left.\frac{\dd}{\dd s}\right|_{s=0}\int_\dom\epsilon \left(1+sH(x)\right)^{\frac{1}{m}}W\dd x
=\int_\dom \epsilon HW\dd x
\end{equation}
is determined by $\Lambda^h_{\epsilon,\gamma}$. From here it is easy to see that we must have
\begin{equation}
\int_\dom \left(\epsilon^{(i)}-\epsilon^{(ii)} \right)UW\dd x=0
\end{equation}
for all $U,W$ such that $\nabla\cdot(\gamma\nabla U)=\nabla\cdot(\gamma\nabla W)=0$. As shown in \cite{SU}, the span of such products $UV$ is dense in $L^2(\dom)$, so this implies that $\epsilon^{(i)}=\epsilon^{(ii)}$.

\paragraph{Acknowledgments:} C.C. was supported by NSTC grant number 112-2115-M-A49-002. The research of T.G. was supported by the Collaborative Research Center, membership no. 1283, Universit\"at Bielefeld. G.N. was supported by Grant-in-Aid for Scientific Research of the Japan Society for the Promotion of Science (19K03554). Work on this project began when C.C. and T.G.  visited the National Center for Theoretical Sciences (NCTS), National Taiwan University. The authors acknowledge the support provided by NCTS during that visit.

\bibliography{pme-2}
\bibliographystyle{plain}

\end{document}